\documentclass[12pt]{article}
\usepackage{amsmath,amssymb,graphicx}

\newtheorem{theorem}{\bf Theorem}[section]

\newtheorem{corollary}[theorem]{\bf Corollary}

\newcommand{\eopf}{$\Box$}

\newenvironment{proof}
{\noindent{\bf Proof.}\ }%
{ \hfill\eopf\par\bigskip}%

\title{Fullerene graphs have exponentially many perfect matchings}
\author{Franti\v sek Kardo\v s$^{a }$\footnote{This author is supported in part by Science and
Technology Assistance Agency under the contract No. APVT-20-004104 and VVGS grant No. I-07-003-31.}, Daniel Kr\' al\!'$^{b }$\footnote{Institute for Theoretical Computer Science ({\sc iti}) is supported by Ministry of Education of the Czech Republic as project $1M0545$.}, Jozef Mi\v skuf\hskip 1pt$^{a }$\footnotemark[1], \\Jean-S\'ebastien Sereni$^{c }$\footnote{This author is supported by the European project {\sc ist fet Aeolus.}}}
\date{}
\begin{document}

\maketitle
\begin{center}
$^a$ Institute of Mathematics, Faculty of Science, University of Pavol Jozef \v{S}af\' arik, Jesenn\'a 5, 041 54
Ko\v{s}ice, Slovakia \\
E-mail: {\tt \{frantisek.kardos,jozef.miskuf\}@upjs.sk}. \\
$^b$ Institute for Theoretical Computer Science ({\sc iti}), Faculty of Mathematics and Physics, Charles University, Malostransk\'e N\'am\v est\'i 25,
118 00 Prague, Czech Republic. \\
E-mail: {\tt kral@kam.mff.cuni.cz}. \\
$^c$ Department of Applied Mathematics ({\sc kam}) and 
{\sc dimatia}, Faculty of Mathematics
and Physics, Charles University, Malostransk\'e N\'am\v est\'i 25,
118 00 Prague, Czech Republic. \\ E-mail: {\tt sereni@kam.mff.cuni.cz}.
\end{center}

\begin{abstract}
A fullerene graph is a planar cubic 3-connected graph with only pentagonal and hexagonal faces. We show that fullerene graphs have exponentially many perfect matchings.
\end{abstract}



\section{Introduction}
Since the discovery of the first fullerene molecule \cite{KHBCS} in 1985, the fullerenes have been objects of interest to scientists all over the world.
The name \emph{fullerenes} was given to cubic carbon molecules in which the atoms are arranged on a sphere in pentagons and hexagons.

Many properties of fullerene molecules can be studied using mathematical tools and results. Thus, \emph{fullerene graphs} were defined as cubic (i.e.~$3$-regular) planar 3-connected graphs with pentagonal and hexagonal faces. Such graphs are suitable models for fullerene molecules: carbon atoms are represented by vertices of the graph, whereas the edges represent bonds between adjacent atoms.

Since all carbon atoms are 4-valent, for every atom precisely one of the three bonds should be doubled. Such a set of double bonds is called a \emph{Kekul\' e structure} in a fullerene. It corresponds to the notion of perfect matchings in fullerene graphs: a \emph{matching} in a graph $G$ is a set of edges of $G$ such that no two edges in $M$ share an end-vertex. A matching $M$ is \emph{perfect} if any vertex of $G$ is incident with an edge of $M$.
Let $M$ be a perfect matching in a fullerene graph $G$. A hexagonal face is \emph{resonant} if it is incident with three edges in $M$. The maximum size of a set of resonant hexagons in $G$ is called the \emph{Clar number} of $G$.


All known general lower bounds for the number of perfect matchings in fullerene graphs are linear in the number of vertices \cite{D1,D2,ZZ}. The best known result asserts that a fullerene graph with $p$ vertices has at least $\left\lceil \frac{3(p+2)}{4}\right\rceil$ different perfect matchings \cite{ZZ}.
On the other hand, the computation of the average number of perfect matchings in fullerene graphs with small number of vertices  \cite{D3} indicates that this number should grow exponentially with $p$.

So far, several special classes of fullerene graphs with exponentially many perfect matchings are known but a general result is missing. Such classes of fullerene graphs either have the special structure of nanotubes \cite{KM}, have high symmetry \cite{D3} or are obtained using specific operations \cite{D4}. 
In this paper, we establish an exponential lower bound on the number of perfect matchings for all fullerene graphs.

\section{Main result}
A \emph{cyclic edge-cut} in a graph $G$ is an edge-cut $E$ such that at least two of the connected components  of $G\setminus E$ contain a cycle. 
A cyclic edge-cut is \emph{trivial}, if one of the components is a cycle. 
No fullerene graph has a cyclic edge-cut of size less than five \cite{D5,KS}.
The fullerene graphs with non-trivial cyclic 5-edge-cuts have a special structure \cite{KS,KM} and the number of perfect matchings in them is known to be exponential \cite{KM}. 
Hence, we focus on fullerene graphs with no non-trivial cyclic 5-edge-cuts. 

\begin{theorem}
Let $G$ be a fullerene graph with $p$ vertices that has no non-trivial cyclic 5-edge-cut. The number of perfect matchings of $G$ is at least $2^{\frac{p-380}{61}}.$
\label{th:1}
\end{theorem}

\begin{proof}
We find a perfect matching $M$ in $G$ such that there are at least $\frac{p-380}{61}$ disjoint resonant hexagonal faces. Since in each such resonant hexagon we can switch the matching to the other edges of the hexagon independently of the other resonant hexagons, the bound will follow immediately.

The dual graph $G^*$ of the graph $G$ is a plane triangulation with 12 vertices of degree 5 and all other vertices of degree 6. Let $U=\{u_1,\dots,u_{12}\}$ be the set of vertices of degree 5. 
Our aim is to construct a set $W$ of vertices of $G^*$ of degree 6 and such that: 
\begin{itemize}
\item the distance between $v$ and $v^\prime$ in $G^*$ is at least $5$ for all $v,v^\prime\in W$, $v\ne v^\prime$;
\item the distance between $v$ and $u$ in $G^*$ is at least $3$ for all $v\in W$ and $u\in U$.
\end{itemize}

We present a greedy algorithm to construct such a set $W$. Initially, we set $W_0=\emptyset$, and we color all the vertices at distance at most 2 from any $u_i$ by the white color. The remaining vertices are colored black. White vertices cannot be chosen as vertices of $W$. For each $u_i\in U$ there are at most 5 vertices at distance one and at most 10 vertices at distance two. Hence, there are at most  $12\cdot (1+5+10)=192$ white vertices initially. 

Until there are some black vertices, we choose a black vertex $v_k$ and add it to the constructed set, i.e.
$W_k:=W_{k-1}\cup\{v_k\}$. We recolor all vertices at distance at most 4 from $v_k$ (including $v_k$) white. Since for any vertex $v$ of degree 6 there are at most $6d$ vertices at distance $d$, there are at most  $1+6+12+18+24=61$ new white vertices. 
This procedure terminates when there are no black vertices. 

Let $W$ be the resulting set $W_k$.
The set $W$ contains at least $\frac{f-192}{61}$ vertices where $f$ is the number of faces of $G$. By Euler's formula, $f=\frac{p}{2}+2$ and thus $|W|\geq\frac{p-380}{122}$.

\begin{figure}[ht]
\centerline{\includegraphics{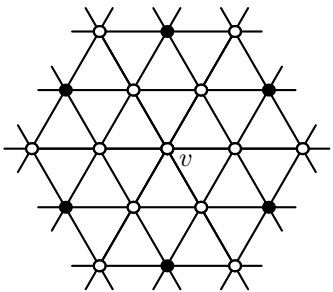}}
\caption{The configuration $R(v)$ and the six vertices in $R^*(v)$.}
\label{fig:conf}
\end{figure}
We next describe how to construct a matching in $G$ with a lot of disjoint resonant hexagons. Given a vertex $v\in W$, let $R(v)$ be the set of vertices at distance at most $2$ from $v$ (see Figure \ref{fig:conf}). The vertices at distance $2$ from $v$ form a cycle of length $12$ in $G^*$. This cycle is  an induced cycle of $G^*$ since $G$ has no non-trivial cyclic $5$-edge-cut. Let $R^*(v)$ be the set formed by the $6$ independent vertices of $R(v)$ drawn with full circles in Figure \ref{fig:conf}. Since $G$ has no non-trivial cyclic $5$-edge-cut, all the vertices in $R^*(v)$ are different and form an independent set in $G^*$.

\begin{figure}[ht]
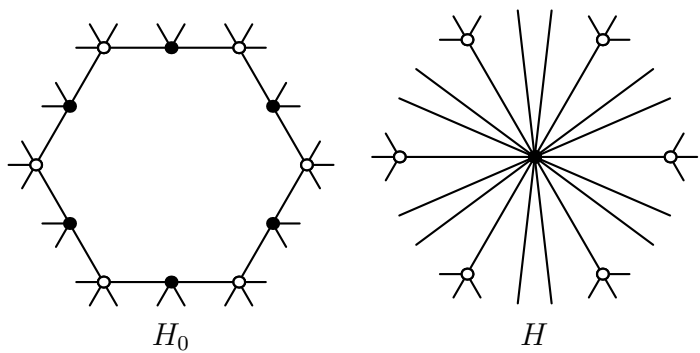

\centerline{
\begin{tabular}{cc}
\includegraphics{conf2.2}&
\includegraphics{conf2.3}\\
$H_0$&$H$
\end{tabular}
}
\caption{The structure of the graphs $H_0$ and $H$.}
\label{fig:h}
\end{figure}
The sets $R^*(v)$ for $v \in W$ are pairwise disjoint since $W$ only contains vertices at distance at least 5. We now modify the graph $G^*$ to planar graphs $H_0$ and $H$. For every vertex $v \in W$, delete $v$ and the six neighbors of $v$. Let $H_0$ be the resulting graph. Further identify the six vertices of $R^*(v)$ (see Figure \ref{fig:h}). The final plane graph is denoted by $H$. 

The Four Color Theorem \cite{4CT11,4CT22,4CT2} asserts the existence of a proper vertex coloring of $H$ using four colors.  
The coloring of $H$ yields a precoloring of $H_0$ such that the six vertices of each set $R^*(v)$ have the same color. Let $c(v)$ be this color.

We extend the precoloring of $H_0$ to a proper coloring of $G^*$. We first color each vertex $v$ by the color $c(v)$. For each $v\in W$, there are only six uncolored vertices inducing a 6-cycle (the vertices adjacent to $v$), and each such vertex has three neighbours colored with $c(v)$ and one vertex colored with a different color. Therefore, for each such uncolored vertex, there are $2$ available colors. Since every cycle of length six is $2$-choosable \cite{ERT,VIZ}, there is an extension of the coloring of $H_0$ to $G^*$.

The $4$-coloring of $G^*$ corresponds to a proper $3$-edge coloring of $G$. To see this, assume that the vertices of the graph $G^*$ are colored with colors $1,2,3$, and $4$. There are edges of 6 different color types: $12$, $13$, $14$, $23$, $24$, and $34$.  Color the edges of $G$ corresponding to the edges of $G^*$ of types $12$ and $34$ (which are pairwise disjoint) by  the color $a$, the edges of $G$ corresponding to the edges of $G^*$ of types $13$ and $24$ by the color $b$, and the remaining edges, i.e. the edges corresponding to the edges of $G^*$ of types $14$ and $23$, by the color $c$. Since the graph $G$ is cubic, each of the color classes $a$, $b$, and $c$ forms a perfect matching of $G$. 

Let $f$ be a face corresponding to a vertex $w$ adjacent to $v\in W$ in $G^*$. Since $w$ has three (non-adjacent) neighbors in $G^*$ colored with the color $c(v)$, the corresponding three non-adjacent edges incident with $f$ are colored with the same color. Hence, the face $f$ is resonant in one of the three matchings formed by the edges of the color $a$, the edges of the color $b$, and the edges of the color $c$. 

There are $6$ such resonant hexagons for the three matchings for each $v\in W$. Altogether, there are $6|W|$ resonant hexagons. Therefore, one of the matchings has at least $2|W|\ge\frac{p-380}{61}$ resonant hexagons.
Observe that the resonant hexagons in one color class are always disjoint: if they were not disjoint, they would correspond to two adjacent neighbors $w$ and $w^\prime$ of some vertex $v\in W$. But the colors assigned to $w$ and $w^\prime$ are different, in particular, the edges corresponding to $vw$ and $vw^\prime$ have different colors. Hence, the hexagons corresponding to $w$ and $w^\prime$ are resonant for different colors $a$, $b$, or $c$. The desired bound on the number of perfect matchings readily follows.
\end{proof}

Theorem \ref{th:1} combined with the bound $15\cdot 2^{\frac{p}{20}-1/2}$ by Kutnar and Maru\v si\v c \cite{KM} on the number of perfect matchings in fullerene graphs with non-trivial cyclic 5-edge cuts yields the following.
\begin{corollary}
Every fullerene graph with $p$ vertices has at least $2^{\frac{p-380}{61}}$ perfect matchings.
\end{corollary}


\begin{thebibliography}{11}
\bibitem{4CT11}  K. Appel and W. Haken, {\em Every planar map is four colorable. Part I. Discharging}, Illinois J. Math. {\bf 21} (1977) 429--490.
\bibitem{4CT22} K. Appel, W. Haken and J. Koch, {\em Every planar map is four colorable. Part II. Reducibility}, Illinois J. Math. {\bf 21} (1977) 491--567.
\bibitem{D1} T.~Do\v sli\'c, {\em On lower bounds of number of perfect matchings in fullerene graphs},  J. Math. Chem. {\bf 24} (1998) 359--364.
\bibitem{D2} T.~Do\v sli\'c, {\em On some structural properties of fullerene graphs},  J. Math. Chem. {\bf 31} (2002) 187--195.
\bibitem{D5} T.~Do\v sli\'c, {\em Cyclical edge-connectivity of fullerene graphs and $(k,6)$-cages}, J. Math. Chem {\bf 33} (2003) 103--112.
\bibitem{D3} T.~Do\v sli\'c, {\em Fullerene graphs with exponentially many perfect matchings}, J. Math. Chem. {\bf 42} (2) (2007) 183--192.
\bibitem{D4} T.~Do\v sli\'c, {\em Leapfrog fullerenes have many perfect matchings}, J. Math. Chem. (2007), doi: 10.1007/s10910-007-9287-x.
\bibitem{ERT} P. Erd\H os, A. L. Rubin, H. Taylor, {\em Choosability in graphs}, Congr. Numer. {\bf 26} (1980) 122--157.
\bibitem{KS} F.~Kardo\v s, R. \v Skrekovski, {\em Cyclic edge-cuts in fullerene graphs}, J. Math. Chem. (2007), doi: 10.1007/s10910-007-9296-9.
\bibitem{KHBCS} H.~W.~Kroto, J.~R.~Heath, S.~C.~O'Brien, R.~F.~Curl and R.~E.~Smalley, {\em C$_{60}$ Buckminsterfullerene}, Nature {\bf 318} (1985) 162--163.
\bibitem{KM} K.~Kutnar, D.~Maru\v si\v c, {\em On cyclic edge-connectivity of fullerenes}, preprint, arXiv:math-CO/0702511.
\bibitem{4CT2} N. Robertson, D. P. Sanders, P. D. Seymour, R. Thomas, {\em The four colour theorem}, J. Combin. Theory Ser. B. {\bf 70} (1997) 2--44.
\bibitem{VIZ} V. G. Vizing, {\em Colouring the vertices of a graph in prescribed colours} (in Russian), Diskret. Anal. {\bf 29} (1976) 3--10.
\bibitem{ZZ} H.~Zhang and F.~Zhang, {\em New lower bound on the number of perfect matchings in fullerene graphs},  J. Math. Chem. {\bf 30} (2001) 343--347.
\end{thebibliography}
\end{document}